\documentclass{jnmp}
\usepackage{amsmath}
\usepackage{graphicx}
\usepackage{latexsym}
\JNMPnumberwithin{equation}{section}
\resetfootnoterule
\allowdisplaybreaks

\newtheorem*{acknowledgments}{Acknowledgments}
\theoremstyle{definition}

\begin{document}

\renewcommand{\evenhead}{Fazal M Mahomed and Asghar Qadir}
\renewcommand{\oddhead}{Linearizability Criteria for a Class of Third Order
Semi-Linear ODEs}
\thispagestyle{empty}
\FirstPageHead{*}{*}{20**}{\pageref{firstpage}--\pageref{lastpage}}{Article}
\copyrightnote{200*}{Fazal M Mahomed and Asghar Qadir}
\Name{Linearizability Criteria for a Class of Third Order
Semi-Linear Ordinary Differential Equations}\label{firstpage}
\Author{Fazal M MAHOMED~$^a$ and Asghar QADIR~$^b$} \Address{$^a$
School of Computational and Applied Mathematics, Centre for
Differential Equations, Continuum Mechanics and Applications\\
University of the Witwatersrand, Wits 2050, South Africa \\~~E-mail: Fazal.Mahomed@wits.ac.za\\[10pt]$^b$ Centre for Advanced Mathematics \& Physics\\
National University of Sciences \& Technology\\
Campus of College of Electrical \& Mechanical Engineering\\
Peshawar Road, Rawalpindi, Pakistan \\~~E-mail: aqadirmath@yahoo.com}
\Date{Received Month *, 200*; Accepted in Revised Form Month *, 200*}
\begin{abstract}
\noindent  Using geometric methods for linearizing
systems of second order cubically semi-linear ordinary differential
equations, we extend to the third order by differentiating the
second order equation. This yields criteria for linearizability of a
class of third order semi-linear ordinary differential equations,
which is distinct from the classes available in the literature. Some
examples are given and discussed.
\end{abstract}
\section{Introduction}
The study of non-linear differential equations initially centred
about procedures to reduce them to linear form. The earliest
attempts tried to approximate the non-linear equation by a linear
one and use the solution to iteratively improve the approximation
[15]. One then has to hope that the series converges. It is never
very clear in this procedure in what regime the essence of
non-linearity will be lost in approximation. Lie used his point
transformations [8] to determine the class of scalar second order
semi-linear ordinary differential equations (odes) that would
transform under it to linear form [9], thus providing {\it exact}
solutions of the odes. Such odes are called {\it linearizable in the
sense of Lie}. More generally, one allows other transformation of
variables. We shall simply call them {\it linearizable}. He
obtained specific criteria for linearizability and obtained the most
general form of the equation to be cubically semi-linear, apart from
obtaining an algebraic classification of the equations that could be
linearized. Some work on classification of systems of two such
equations was done [17] and it was found that there are five classes
of such linearizable systems. Some further work proved that for the
scalar third order there are three classes [10]. However, till recently
nothing had been obtained about explicit linearizability criteria of
systems of second order, or of third order odes.

Chern [2,3] used contact transformations to discuss the
linearizability of scalar third order odes reducible to the linear forms $u'''=0$
and $u'''+u=0$. Point transformations were used by Grebot [5,6] but
were restricted to the class of transformations $t=\phi(x)$,
$u=\psi(t,x)$. This work was generalized by Neut and Petitot, and
Ibragimov and Meleshko [16,7] to the third class known to exist
[17], $u'''+\alpha (t)u=0$. They obtained a larger class of third
order semi-linear odes that is linearizable. Meleshko [14] also
considered arbitrary point transformations to reduce third order
odes of the form $y'''=F[y,y',y'']$ and used the $\partial/\partial
t$ symmetry to reduce the order to two. He then used the Lie
linearizability criteria to determine the linearizability of the
scalar third order equation.

Using the connection between the isometry algebra and the symmetries
of the system of geodesic equations [4], linearizability criteria
were stated for a system of second order quadratically semi-linear
odes, of a class that could be regarded as a system of geodesic
equations [11]. The criteria came from requiring that the
coefficients in the equations, regarded as Christoffel symbols,
yield a zero curvature tensor. We use the projection procedure of
Aminova and Aminov [1], which appeals to the fact that the geodesic
equations do not depend on the geodetic parameter, to reduce the
dimension by one and convert from a quadratically semi-linear system
to a more general cubically non-linear system [12]. When applied to
a system of two dimensions we get a single cubically semi-linear
ode. The linearizability criteria so obtained are exactly the Lie
criteria! Applied to a system of three dimensions, one obtains a
system of two cubically semi-linear odes with extended Lie criteria.
The system of odes so obtained is in the class of maximally
symmetric equations (out of the total of five classes mentioned
earlier). Some other classes come from projection of the original
space to lower dimensional spaces.

In this paper we consider the linearizability criteria of a class of scalar third order odes obtained by
differentiating the cubically semi-linear system of second order odes and stating
criteria of linearizability for the class so obtained. (It is worth
mentioning that we obtain the same criteria if we first
differentiate the system and then project it.) The class obtained
here is distinct from all the earlier classes considered as there is
no guarantee that there will be three arbitrary constants appearing
in the solution.

The plan of the paper is as follows. In the next section we present
the geometrical notation used and briefly review the essential
results used in the sequel. In section 3 we derive the
linearizability criteria for the scalar third order ode in the form involving the
second derivative and in a form quintic in the first derivative. In
the following section we present some illustrative examples.
Finally, in section 5 we give a brief summary and discussion of the
results.
\section{Notation and review}
The following are well-known and can be found in text books. We use
the Einstein summation convention that repeated indices are summed
over the entire range of the index. Thus, $A^aB_a$ stands for
$\Sigma_{a=1}^n A^aB_a$. The metric tensor will be represented by
the symmetric (non-singular) matrix $g_{ij}$ and its inverse by
$g^{ij}$. The Christoffel symbols are given by
$$
\Gamma^i_{jk}=\frac12 g^{im}(g_{jm,k}+g_{km,j}-g_{jk,m}),\eqno{(1)}
$$
where $,k$ stands for partial derivative relative to $x^k$, etc. The
Christoffel symbols are symmetric in the lower pair of indices,
$\Gamma^i_{jk}=\Gamma^i_{kj}$.

In this notation, the system of $n$ geodesic equations is
$$
\ddot x^i+\Gamma^i_{jk}\dot x^j\dot x^k=0,\quad i,j,k=1,\ldots,n,\eqno{(2)}
$$
where $\dot x^i$ is the derivative relative to the arc length
parameter $s$ defined by $ds^2=g_{ij}dx^i dx^j$.

The Riemann tensor is defined by
$$
R^{i}_{jkl}=\Gamma^i_{jl,k}-\Gamma^i_{jk,l}+\Gamma^i_{mk}\Gamma^m_{jl}-\Gamma^i_{ml}
\Gamma^m_{jk}, \eqno{(3)}
$$
and has the properties
$$
R^i_{jkl}=-R^i_{jlk},\eqno{(4)}
$$
$$
R^i_{jkl}+R^i_{klj}+R^i_{ljk}=0, \eqno{(5)}
$$
and
$$
R^i_{jkl;m}+R^i_{jlm;k}+R^i_{jmk;l}=0. \eqno{(6)}
$$
The Riemann tensor in fully covariant form is
$$
R_{ijkl}=g_{im}R^m_{jkl},\eqno{(7)}
$$
and satisfies
$$
R_{ijkl}=-R_{jikl}.\eqno{(8)}
$$
The linearization criteria of [11] are $R_{ijkl}=0$.

For a general system of two odes for two functions, if it can be
regarded as a system of geodesic equations,
$$
x''=a(x,y)x'^2+2b(x,y)x'y'+c(x,y)y'^2,
$$
$$
y''=d(x,y)x'^2+2e(x,y)x'y'+f(x,y)y'^2,\eqno{(9)}
$$
the Christoffel symbols in terms of these coefficients are
$$
\Gamma_{11}^1=-a, \Gamma^1_{12}=-b, \Gamma_{22}^1=-c,
$$
$$
\Gamma_{11}^2=-d, \Gamma^2_{12}=-e, \Gamma_{22}^2=-f.\eqno{(10)}
$$
The linearizability criteria are
$$
a_y-b_x + be -cd = 0 , b_y -c_x+ (ac - b^2) + (bf - ce) = 0 ,
$$
$$
d_y-e_x - (ae - bd) - (df - e^2) = 0 , (b + f)_x = (a + e)_y , \eqno{(11)}
$$
with constraints on the metric coefficients: $g_{11}=p , g_{12}=g_{21}=q , g_{22}=r$,
$$
p_x = -2(ap + dq) , q_x = -bp - (a + e)q - dr , r_x = -2(bq + er) ,
$$
$$
p_y = -2(bp + eq) , q_y = -cp - (b + f)q - er , r_y = -2(cq + fr) . \eqno{(12)}
$$
The compatibility of this set of six equations gives the above four
linearization conditions (11). For a system of three equations we
get eighteen such equations out. In general there are $n^2(n+1)/2$,
$n\ge 2$. One now obtains the required linearizing transformation by
regarding the variables in which the equation is linearized as
Cartesian, thus having $g_{11}=g_{22}=1$ and $g_{12}=g_{21}=0$, and
looking for the coordinate transformation yielding the original
metric coefficients [11].

Following Aminova and Aminov [1] and projecting the system down by
one dimension, the geodesic equations become [12]
$$
{\ddot x}^a+A_{bc}{\dot x}^a{\dot x}^b{\dot x}^c+B^a_{bc}{\dot
x}^b{\dot x}^c+C^a_b{\dot x}^b+D^a=0,\quad a=2,\ldots,n,\eqno{(13)}
$$
where the dot now denotes differentiation with respect to the
parameter $x^1$ and the coefficients in terms of the
$\Gamma^a_{bc}$s are
$$
A_{bc}=-\Gamma^1_{bc}, B^a_{bc}=\Gamma^a_{bc}-2\delta^a_{(c}\Gamma^1_{b)1},
$$
$$
C^a_b=2\Gamma^a_{1b}-\delta^a_b\Gamma^1_{11}, D^a=\Gamma^a_{11},
a,b,c=2,\ldots,n,\eqno{(14)}
$$
which is cubically semi-linear. For $n=2$, writing $h=a-2e, g=f-2b$,
we get the Lie criteria for the scalar equation involving the two
auxiliary variables $b$ and $e$, the former of which corresponds to
our $b$ and the latter to our $-e$:
$$
3b_x-2g_x+h_y-3be-3cd=0, b_y-c_x+b^2+bg+ch-ce=0,
$$
$$
d_y+e_x-bd-dg-e^2+eh=0, 3e_y-2h_y-g_x+3be+3cd=0, \eqno{(15a)}
$$
which yield, in terms of the four coefficients of the cubically
semi-linear equation [12],
$$
3(ch)_x+3dc_y-2gg_x-gh_y-3c_{xx}-2g_{xy}-h_{yy}=0,
$$
$$
3(dg)_y+3cd_x-2hh_y-hg_x-3d_{yy}-2h_{xy}-g_{xx}=0. \eqno{(15b)}
$$
We will use these criteria to determine the linearizability of the
scalar third order semi-linear ode.

\section{Derivation of the linearizability criteria for the third order
ode}
\quad
On differentiating (13) for the single differential equation,
writing the independent variable as $x$ and the dependent variable
as $y$, we get
$$
y'''+(3cy'^2-2gy'+h)y''+c_yy'^4+(c_x-g_y)y'^3-(g_x-h_y)y'^2+(h_x-d_y)y'-d_x=0,
\eqno{(16)}
$$
which is a total derivative. We could retain the equation in this
form or use the original to remove the second order term and write
the equation as quintically non-linear in the first derivative:
$$
y'''-3c^2y'^5+(5cg+c_y)y'^4-(4ch+2g^2+g_y-c_x)y'^3
$$
$$
-(2gd+h^2-h_x+d_y)y'^2-(2dg+h^2+d_y-h_x)y'+(dh-d_x)=0. \eqno{(17)}
$$
As we show shortly, for the linearizability criteria it makes no
difference which form we use. However the second, quintic, form is
no longer a total derivative of the original equation except in the
linear case. For purposes of comparison with the forms in the
literature one needs (16) but for comparison with the previous work
on the cubically semi-linear systems of odes (17) may be more
convenient.

For (17) the general form of the equation is
$$
y'''-\alpha y'^5+\beta y'^4-\gamma y'^3+\delta y'^2-\epsilon y'+\phi=0. \eqno{(18)}
$$
In this case for compatibility of (17) and (18) the conditions are
$$
\alpha=3c^2, \eqno{(19)}
$$
$$
\beta=5cg+c_y, \eqno{(20)}
$$
$$
\gamma=4ch+2g^2+g_y-c_x, \eqno{(21)}
$$
$$
\delta=3cd+3gh+h_y-g_x, \eqno{(22)}
$$
$$
\epsilon=2dg+h^2+d_y-h_x, \eqno{(23)}
$$
$$
\phi=dh-d_x. \eqno{(24)}
$$
On inversion, the first four equations provide definitions for the
four coefficients in (13), (14) and (15b) for $n=2$, and the next
two provide consistency constraint conditions:
$$
c={\sqrt \alpha /\sqrt 3}, \eqno{(25)}
$$
$$
g=(\beta -c_y)/5c, (c \ne 0), \eqno{(26)}
$$
$$
h=(\gamma -2g^2-g_y+c_x)/4c, (c \ne 0), \eqno{(27)}
$$
$$
d=(\delta -3gh-h_y+g_x)/3c, (c \ne 0), \eqno{(28)}
$$
$$
\epsilon=2dg+h^2+d_y-h_x, (c \ne 0), \eqno{(29)}
$$
$$
\phi=dh-d_x, (c \ne 0). \eqno{(30)}
$$
In the case $c=0$, clearly $\alpha=\beta=0$, and the equation
becomes cubically non-linear. Now there can be different choices of
$g$ for a given $\gamma$, $h$ for a given $\delta$ and choice of
$g$, and $d$ for a given $\epsilon$ and choices of $g$ and $h$.
There is then only one consistency condition for the various
choices. An example will be given in the next section to illustrate
this case.

Hence we have the following theorem. \\
{\bf Theorem 1:} {\it Equation (18) is linearizable if the
linearizability criteria (15) are satisfied, where the coefficients
are given by (25) - (30) ($c \ne 0$), with $h=a-2e, g=f-2b$}.

The general form corresponding to (16) is
$$
y'''+(A_2y'^2-A_1y'+A_0)y''+B_4y'^4-B_3y'^3+B_2y'^2-B_1y'+B_0=0. \eqno{(31)}
$$
The identification of most of the coefficients is {\it easier} in
this form. We have
$$
c=A_2/3, g=A_1/2, h=A_0, \eqno{(32)}
$$
but $d$ is not obtained from here. The constraint conditions arising
are:
$$
B_4=A_{2y}/3, B_3=A_{1y}/2-A_{1x}/3, B_2=A_{0y}-A_{1x}/2, \eqno{(33)}
$$
which are also easier to check than the corresponding equations in
the other form. There are two differential equations for $d$ that
yield its value up to an arbitrary constant
$$
d=-\int B_2dx+k(y)=\int (B_1-A_{0x})dy+l(x). \eqno{(34)}
$$
The constant can be determined by requiring consistency of the
Christoffel symbol $\Gamma^2_{11}$. {\it This} is why it becomes
more difficult to compute with this form of the equation. Thus we
have the following theorem. \\
{\bf Theorem 2:} {\it Equation (31) is linearizable if the
linearizability criteria (15) are satisfied, where the coefficients
are given by (32) - (34) ($c \ne 0$), with $h=a-2e, g=f-2b$, after
requiring consistency of the Christoffel symbols with the deduced
metric coefficients.}

The four linearizability conditions (15) are stated in terms of the
6 coefficients $a,...,f$ and not the four coefficients $c,d,g,h$.
Thus there is degeneracy in the choices available. Any choices of
$a$ and $e$ for a given $h$, or $f$ and $b$ for a given $g$,
compatible with the metric coefficient relations (12) are
permissible. For each such choice we would get corresponding
linearizability conditions.

For example, assume that $b=e=0$, then
$$
p_x=-2(ap+dq), p_y=0, q_x=-(aq+dr), \eqno{(35)}
$$
$$
q_y=-(cp+fq), r_x=0, r_y=-2(cq+fr). \eqno{(36)}
$$
This is a consistent set of requirements. The other four conditions
become
$$
h_y=cd, c_x=ch, d_y=dg, g_x=h_y. \eqno{(37)}
$$
The first two of these are
$$
4h_y=\delta-3gh+g_x,
$$
$$
\alpha_x/2\alpha=[50\sqrt 3\alpha^2\gamma + 25\alpha^{3/2}+
\sqrt 3(-2\sqrt 3\beta +\alpha_y)^2
$$
$$
+5\sqrt 3(-\alpha_y^2+\alpha \alpha_{yy}+2\sqrt{3\alpha}(\beta
\alpha_y-\alpha \beta_y))]200\alpha^{5/2}. \eqno{(38)}
$$
Hence the first and last of the four conditions are
$$
g_x=h_y=\delta/3-gh. \eqno{(39)}
$$
Since $h$ itself depends on the second derivative of $\alpha$
relative to $y$, the last condition above is a third order
derivative constraint.

Alternatively, choosing $a=f=0$, we would get
$$
p_x=-2dq, p_y=-2(bp+eq), q_x=-(bp+eq+dr),
$$
$$
q_y=-(cp+bq+er), r_x=-2(bq+er), r_y=-2cq, \eqno{(40)}
$$
and
$$
g_x=h_y=2cd-gh/2, 2c_x+g_y=-2(b^2+ce), 2d_x+g_x=-(dg+h^2/2).
\eqno{(41)}
$$
These are clearly more complicated and will not be used.

\section{Examples}
We present some examples of third order equations that can be
linearized by our procedure.

{\bf 1.} Here we choose $e=f=0$ so that $g=-2b, h=a$. Further, let
us try to obtain a lower degree of the equation in the form that
there is no second derivative. As pointed out earlier this requires
that $c=0$. To construct this example we choose the metric
coefficients $p=A(x)y^2, q=B'(x)y, r=2B(x)$. Writing $pr-q^2=\Delta
y^2$ and $A^2+A'B'/2-AB''=\Lambda$, we get $h=a=-(ln \sqrt \Delta)',
g=-2b=2/y, c=0, d=\Lambda y/\Delta, e=f=0$, where the prime here
refers to differentiation with respect to only $x$ and {\it not}
relative to both $x$ and $y$. Notice that the degeneracy of the case
when $c=0$ has been taken care of by selecting the metric
coefficients appropriately. Here, choosing $A=c_1e^{-kx},
B=c_2e^{-kx}$ the equation becomes
$$
y'''-6y'^3/y^2+8ky'^2/y-(k^2-5l)y'+kly=0, \eqno{(42)}
$$
where $c_1=2lc_2$. This is amenable to reduction by the method of
[14] and is not of a form given in [16] or [7].

{\bf 2.} Here we construct the example by differentiating the
linearizable equation
$$
y''+xy'^3+\frac{2}{x}y'=0, \eqno{(43)}
$$
which yields
$$
y'''-3x^2y'^5-7y'^3-\frac{6}{x^2}y'=0. \eqno{(44)}
$$
This has $c = x, d = 0, g = f - 2b = 0$ and $h = a - 2e = 2/x$. We
choose $b = 0$ and $e = -1/x$. Then $f = 0$ and $a = 0$. It is
easily verified that the conditions are met and that it is not in
any of the classes of [16], [7] or [14] and cannot be linearized
according to the methods therein. It can easily be reduced to the
simplest linear form by the transformation $u = x \cos y, v = x \sin
y$. Here the solution is in terms of two arbitrary constants.

{\bf 3.} In this example the equation taken is
$$
y'''-3x^2y'^5/y^4-3xy'^4/y^3+6y'^3/y^2+6y'^2/xy-6y'/x^2=0.
\eqno{(45)}
$$
Since the equation is non-trivial in that it is not in the class of
Neut and Petiot [16], Ibragimov and Meleshko [7] or Meleshko [14],
is not a total derivative, nor can it be solved by other simple
means, we give the steps in some detail. From (25) - (30), we see
that $c=-x/y^2, d=0, h=2/x, g=1/y$. Choosing $a=b=0$ we have $f=g,
e=-h/2$. The metric coefficients are then given by (12) to be
$$
p=y^2(1+y^{-4})/2, q=xy(1-y^{-4})/2, r=x^2(1+y^{-4})/2. \eqno{(46)}
$$
Writing the Cartesian coordinates as $(u,v)$, the coordinate
transformations are given by ([11], eqs. (32) - (34))
$$
u_x^2+v_x^2=y^2+y^{-2}, u_xu_y+v_xv_y=x(y-y^{-3}),
u_y^2+v_y^2=x^2(y^2+y^{-2}). \eqno{(47)}
$$
This set of equations is easily solved by setting
$$
u_x=y, u_y=x, v_x=y^{-1}, v_y=-xy^{-2}, \eqno{(48)}
$$
which yields
$$
u=xy, v=xy^{-1},  \eqno{(49)}
$$
which are the linearizing transformations. Thus the solution is
$$
Axy+Bx/y=1, \eqno{(50)}
$$
where $A$ and $B$ are arbitrary constant real numbers. Notice that
this equation is quintic in the first derivative and does depend on
$x, y, y'$, which is not in the classes of [16], [7] or [14].

\section{Concluding remarks}
\quad
We have shown that the criteria for linearizability of a third order
semi-linear ode that were available in the literature [16, 7]
exclude classes of third order equations that can be obtained by
differentiating second order linearizable equations. By
differentiating the general second order cubically semi-linear ode
we can check the linearizability of the class of third order
semi-linear equations that are products of a quadratic factor in the
first order times the second order term and quartic in the first
order terms. The class can be more conveniently expressed as a
quintic in the first order. In this case it is no longer a total
derivative of the original second-order equation (except in the
linear case) so that it is now hidden that this third order equation
arises from a second order one. As such, our procedure gives a {\it
non-classical solution}. The criteria are explicit for the quintic
form but become degenerate in the case that the factor reduces to a
cubic. It is worth noting that this procedure does not allow for a
quartic factor only. Some examples were given.

The procedure presented here leaves freedom of choice of the
coefficients of the system of quadratic equations that yields the
cubic. The reason is that there were six coefficients of the
two-dimensional quadratic system, while only four combinations of
them enter into the third order equation. There are, indeed, six
coefficients in the general third order equation of the desired
type. The extra two coefficients yield additional constraint
equations for the compatibility of the cubic form with the
derivative of the quadratic form. We then only have to check that
the Lie criteria are satisfied by the corresponding second order
equation.

Since [16] and [7] used the classic Lie procedure and obtained all
the classes possible for the third order ode to be linearizable by
point transformations, the question arises how we have got
additional classes. The answer is that we are not using point
transformations {\it directly}. We use point transformations and
then differentiate. This could be something like contact
transformations (as the derivative is also involved) but it is not
quite a contact transformation procedure either. Note that our
procedure does not guarantee that there will be three arbitrary
constants. We only guarantee two. The procedure of [16] does
guarantee three. It would be interesting if a procedure could be
found that involves only a single constant.

Our procedure is also interesting in that it can yield criteria for
systems of third order odes. Work is in progress for a system of
third order cubically semi-linear system [13]. The equivalent of the
present work for a system is in principle possible but appears too
messy at present. Better algebraic computing software may make it
more manageable.

Another line of work that will be useful is to find methods for the
degeneracy in our procedure to be removed. Alternatively, we need to
determine all classes that allow for the degeneracy, where it cannot
be removed.
\begin{acknowledgments}
AQ is most grateful to DECMA and the School of Computational and
Applied Mathematics, University of the Witwatersrand and for some
useful comments by Profs. P.G.L. Leach, S. Meleshko and R. Popovych.
\end{acknowledgments}%

\label{lastpage}

\begin{thebibliography}{99}\small
\bibitem{A-A} \textsc{Aminova A V} and  \textsc{Aminov N A M}, Projective geometry of
systems of differential equations: general conceptions, {\it Tensor N S}
{\bf 62} (2000), 65--86.
\bibitem{C-2} \textsc{Chern S S}, The geometry of the differential equation
$y'''=F(x,y,y,y'')$, Sci, Rep. Nat. Tsing Hua Univ. 4 (1940), 97--111.
\bibitem{C-1} \textsc{Chern S S}, Sur la geometrie d'une equation differentielle du
troiseme orde, {\it C.R. Acad. Sci. Paris}, (1937) 1227.
\bibitem{F-M-Q} \textsc{Feroze T}, \textsc{Mahomed F M} and \textsc{Qadir A}, The connection between
isometries and symmetries of geodesic equations of the underlying
spaces,  {\it Nonlinear Dynamics} {\bf 45} (2006), 65.
\bibitem{G-1} \textsc{Grebot G}, The linearization of third order ODEs, preprint 1996.
\bibitem{G-2} \textsc{Grebot G}, The characterization of third order ordinary
differential equations admitting a transitive fibre-preserving point
symmetry group, {\it J. Math. Anal. Applic.} {\bf 206} (1997), 364--388.
\bibitem{I-M} \textsc{Ibragimov N H} and \textsc{Meleshko S V}, Linearization of third-order
ordinary differential equations by point and contact transformations
{\it J. Math. Anal. Applic.} {\bf 308} (2005), 266--289.
\bibitem{L-} \textsc{Lie S} {Theorie der Transformationsgruppen} {\it Math. Ann.} {\bf 16} (1880), 441.
\bibitem{L-1} \textsc{Lie S} Klassifikation und Integration von gew\"onlichen
Differentialgleichungenzwischen $x$, $y$, die eine Gruppe von
Transformationen gestaten. {\it  Arch. Math.} {\bf VIII, IX} (1883), 187.
\bibitem{M-L} \textsc{Mahomed F M} and \textsc{Leach P G L} Symmetry Lie Algebras of $n$th Order
Ordinary Differential Equations {\it J. Math Anal Applic} {\bf151} (1990), 80.
\bibitem{M-Q} \textsc{Mahomed F M} and \textsc{Qadir A}, Linearization criteria for a
system of second-order quadratically semi-linear ordinary
differential equations, {\it Nonlinear Dynamics} {\bf 48} (2007), 417.
\bibitem{M-Q} \textsc{Mahomed F M} and \textsc{Qadir A}, Linearization Criteria for Systems of
Cubically Semi-Linear Second-Order Ordinary Differential Equations, preprint 2007.
\bibitem{M-N-Q} \textsc{Mahomed F M}, \textsc{Naeem I} and \textsc{Qadir A}, Conditional linearizability criteria
for a system of third-order ordinary differential equations, preprint 2007.
\bibitem{M-} \textsc{Meleshko S V}, On linearization of third-order ordinary
differential equations {\it J. Phys. A.: Math. Gen. Math.} {\bf 39} (2006), 15135--45.
\bibitem{N-} \textsc{Nayfeh A H}, `Perturbation Methods, Wiley-Interscience, New York, 1973.
\bibitem{N-P} \textsc{Neut S} and \textsc{Petitot M}, La g\'eom\'etrie de
l'\'equation $y'''=f(x,y,y',y'')$, {\it C.R. Acad. Sci. Paris S\'er
I} {\bf 335} (2002), 515--518.
\bibitem{W-M} \textsc{Wafo Soh C} and \textsc{Mahomed F M}, Symmetry breaking for a
system of two linear second-order ordinary differential equations,
{\it Nonlinear Dynamics}  {\bf22} (2000), 121.
\end{thebibliography}
\end{document}